\documentclass{IOS-Book-Article}     
\usepackage{amssymb} 
\usepackage{graphicx}

\usepackage[dvips]{rotating} 

\usepackage[msc-links]{amsrefs}

\usepackage{cite}

\usepackage{hyperref}

\hypersetup{
  colorlinks   = true, 
  urlcolor     = blue, 
  linkcolor    = blue, 
  citecolor   = red 
}

\input xy   \xyoption{all}

\newtheorem{thm}{Theorem}
\newtheorem{prop-theorem}{Proposed Theorem}
\newtheorem{prop-alg}{Proposed Algorithm}

\newtheorem{lem}{Lemma}

\newtheorem{exa}{Example}

\newtheorem{prob}{Problem}

\def\a{{\alpha }} 
\def\C{\mathbb C}
\def\Q{\mathbb Q}
\def\R{\mathbb R}

\def\M{\mathcal M}
\def\X{\mathcal X}
\def\H{\mathcal H}

\newcommand\Aut{\mbox{Aut }}
\def\p{\mathfrak p}

\def\mH{\mathcal H}

\def\P{\mathbb P}

\begin{document}
\begin{frontmatter}          
%
\title{The case for  superelliptic curves}
\runningtitle{Superelliptic curves}

\author[A]{\fnms{L.} \snm{Beshaj}} 
\author[B]{\fnms{T.} \snm{Shaska}}
\author[C]{\fnms{E.} \snm{Zhupa}} 

\runningauthor{Beshaj/Shaska/Zhupa}
\address[A]{Department of Mathematics, \\ Oakland University, Rochester, MI, USA; \\ E-mail: beshaj@oakland.edu}
\address[B]{Department of Mathematics, \\ Oakland University, Rochester, MI, USA; \\ E-mail: shaska@oakland.edu}
\address[C]{Department of Computer Science, \\ University of Information Science and Technology "St. Paul the Apostle", \\ Ohrid, Republic of Macedonia;  \\  E-mail: eustrat.zhupa@uist.edu.mk}

\begin{abstract}
There is a natural question to ask whether the rich mathematical theory of the hyperelliptic curves can be extended to all superelliptic curves.  Moreover, one wonders if all of the applications of hyperelliptic curves such as  cryptography, mathematical physics, quantum computation, diophantine geometry, etc can carry over to the superelliptic curves.  
In this short paper we make the case that the superelliptic curves are exactly the curves that one should study 
\end{abstract}


\end{frontmatter}


\section{Introduction}


This lecture is intended to be a motivation for the study of superelliptic curves which were the main focus of the NATO Advanced Study Institute held in Ohrid, Macedonia in the Summer 2014.  While the scope of interesting mathematical problems related to superelliptic curves is very broad and the applications include different areas of sciences, in this paper we focus on few of the arithmetic problems related to the moduli space of curves, automorphisms groups, minimal models of curves, and rational points on curves with the intention to emphasize the methods of extending the knowledge of hyperelliptic curves to all superelliptic curves. 

A superelliptic curve, simply stated,  is an algebraic curve with equation $y^n = f(x)$ defined over some field $k$.  The well known hyperelliptic curves, $y^2=f(x)$,  are the simplest case of superelliptic curves.  There is a natural question to ask whether the rich mathematical theory of the hyperelliptic curves can be extended to all superelliptic curves.  Moreover, one wonders if all of the applications of hyperelliptic curves such as  cryptography, mathematical physics, quantum computation, diophantine geometry, etc can carry over to the superelliptic curves.  

In this short paper we try to make the point that  superelliptic curves are the exactly the right family of curves to study for several reasons, among them: \\

i)  they are the majority of curves with nontrivial automorphism group, for a fixed genus $g > 2$. 

ii) we can fully determine their automorphism group and write down their explicit equations

iii) they are natural generalizations of the hyperelliptic curves which are well studied. \\

In section 2 we give a very brief description of the terminology and notation that will be used here.  We assume that the reader is familiar with most of the basic theory  of algebraic curves.

A generic curve of genus $g > 2$ has no nontrivial automorphisms.  Curves with non-trivial automorphisms   constitute the singular locus of the moduli space $\M_g$, for $g \geq 3$.  This singular locus   contains several components which are determined by the automorphism group of the curves and the way it acts on the curve (the signature).   The majority of these loci in $\M_g$  correspond to superelliptic curves. In section 2  we show exactly what happens  for $g=4$, where about 70-80 \% of all curves with non-trivial automorphisms are superelliptic curves. It would be interesting to investigate what  happens to this ratio when $g \to \infty $.  Similar to what has been done in \cite{serdica} for hyperelliptic curves, one can attempt to do for every "root" case for the superelliptic curves.  In the case of genus $g=4$ this would correspond to cases 36, 38, 39, and 40 of the diagram.

Superelliptic curves are interesting from another viewpoint.  They are the only curves for which we can write down equations when we know the automorphism group of the curve.  Hence, they are the "only" curves for which we can have a systematic study of their arithmetic properties or even consider arithmetic questions.   We illustrate the genus $g=4$ with Table~\ref{tab_eq}. 

The following problems are some of the main problems that connect the theory of the moduli spaces $\M_g$ with the arithmetic properties of the algebraic curves.  At the current time they can be asked only for superelliptic curves. 

\begin{prob}
Given a superelliptic curve $C$ defined over $\C$, determine an equation for $C$ over its minimal field of definition; see \cite{nato-8} for details.   
\end{prob}

The first problem is studied for small $g$ in the hyperelliptic case by Mestre  and for curves with extra automorphisms in \cite{bi-3, sh-v} and many other references. There are no explicit works for superelliptic curves (non-hyperelliptic), but theoretically things should work the similarly as in the hyperelliptic case.

\begin{prob}
Given a superelliptic curve $C$ defined over an algebraic number field $K$, find an equation for $C$ with 

i) minimal height; see \cite{nato-8} for details.  

ii)  minimal absolute height; see \cite{nato-8} for details. 
\end{prob}

This problem has been studied by Beshaj in \cite{nato-6} and \cite{nato-8} where the notion of minimal height and minimal absolute height are defined. 

\begin{prob}
Create a database of all superelliptic curves defined over $\Q$ for small genus and small height.  
\end{prob}

This problem and some of it ramifications have been suggested in \cite{zhupa} where it is provided a computer program how to write down all the superelliptic curves and their parametric equations for a fixed genus $g>2$. 

It is reasonable for superelliptic  curves   to study their arithmetic properties or to count them according to the moduli height as suggested in \cite{nato-8}.





\def\H{\mathcal H}
\newcommand\A{\mathcal A}
\newcommand\B{\mathcal B}
\def\L{\mathcal L}
\newcommand\E{\mathcal E}

\newcommand\Kc{\mathcal K}

\def\Y{\mathfrak Y}
\def\v{\mathfrak v}
\def\u{\mathfrak u}
\def\w{\mathfrak w}

\newcommand\m{\mathfrak m}
\newcommand\an{\mathfrak a}
\newcommand\bn{\mathfrak b}
\newcommand\en{\mathfrak e}
\newcommand\hn{\mathfrak h}
\newcommand\n{\mathfrak n}
\newcommand\HS{\mathfrak H}

\newcommand\pf{\mathfrak p}

\newcommand\T{\theta}

\newcommand\embd{\hookrightarrow}

\newcommand\G{\bar{G}}
\newcommand\bG{\bar G}

\newcommand\D{\Delta}
\newcommand\e{\varepsilon}
\newcommand\normal{\triangleleft }
\newcommand\car{\overset c \normal}
\newcommand\emb{\hookrightarrow }

\def\a{\alpha}
\def\t{\tau}
\newcommand\s{\sigma}
\newcommand\st{\star}
\def\d{{\delta }}

\newcommand\om{\omega}
\newcommand\vT{\vartheta}

\newcommand\ep{\epsilon}
\newcommand\xs{{\Bbb o}}

\newcommand\sem{\rtimes}

\def\Pic{\mbox{Pic }}
\def\Jac{\mbox{Jac }}
\def\mod{\mbox{ mod }}
\def\deg{\mbox{deg }}
\def\Sel{\mbox{Sel}}
\def\mH{\mathcal H}
\def\Ha{\mbox{H^{\mathbb A}}}

\def\embd{\hookrightarrow}

 
\def\diag{\mbox{diag }}
\def\ch{\mbox{char }}
\def\sem{{\rtimes}}

\def\p{\mathfrak p}
\def\O{\mathcal O}

\def\div{\mbox{div}}
\def\ord{\mbox{ord}}
\def\zz{\zeta}
\def\z{\omega}

\def\Z{\mathbb Z}
\def\bZ{\mathbb Z}
\def\Q{\mathbb Q}
\def\C{\mathbb C}
\def\bP{\mathbb P}
\def\cB{\mathcal B}
\def\cA{\mathcal A}
\def\L{\mathcal L}
\def\cR{\mathcal R}
\def\H{\mathcal H}
\def\M{\mathcal M}
\def\N{\mathcal N}
\def\w{\widetilde}
\def\l{\lambda}
\def\s{\sigma}
\def\a{\alpha}
\def\b{\beta}
\def\p{\mathfrak p}
\def\P{\mathcal P}
\def\e{\varepsilon}
\def\iso{\equiv}
\def\sem{{\rtimes}}

\def\bG{\overline G}
\def\g{\gamma}
\def\bg{\bar \gamma}
\def\u{\mathfrak u}
\def\k{\bar k}
\def\iso{{\, \cong\, }}
\def\nor{{\,  \vartriangleleft \, }}
\def\<{\langle}
\def\>{\rangle}
\def\emb{\hookrightarrow }
\def\rank{\mbox{rank }}

\def\Y{\mathcal Y}
\def\z{\omega}
\def\normal{\triangleleft}
\def\D{\Delta}

\section{Preliminaries on  curves}

By a \textit{curve} we mean a complete reduced algebraic curve over $\C$ which might be singular or reducible.  A \textit{smooth curve} is implicitly assumed to be irreducible. The basic invariant of a smooth curve $C$ is its genus which is half of the first Betti number of the underlying topological space.  We will denote the genus of $C$ by $g(C) = \frac 1 2 \, \rank ( H^1 (C, \Z) )$.

Let $f: \X \to \Y$ be a non-constant holomorphic map between smooth curves $\X$ and $\Y$ of genera $g$ and $g^\prime$.  For any $q\in \X$ and $p=f(q)$ in $\Y$  chose local coordinates $z$ and $w$ centered at $q$ and $p$ such that $f$ has the standard form $w= z^{\nu (q)}$. Then, for any $p$ on $C^\prime$ define \[ f^\star (p) = \sum_{q \in f^{-1} (p) } \nu ( q ) q.\]
If $D$ is any divisor on $\X$ then define $f^\star (D)$ to be the divisor on $\X$ by extending the above $f^\star (p)$. 
The degree of $n$ the divisor $f^\star (p)$ is independent of $p$ and is called the \textit{degree} of the map $f$.  The \textit{ ramification divisor } $R$ on $C$ of the map $f$ is defined by
\[ R = \sum_{q\in C} \left( \nu (p) -1 \right) q \]
The integer $\nu (q) -1$ is called the \textit{ramification index} of $f$ at $q$.   For any meromorphic differential $\phi$ on $C^\prime$ we have \[ \left( f^\star (\phi) \right) = f^\star \left( ( \phi) \right) + R\] Counting degrees we get the Riemann-Hurwitz formula 
\[ 2g-2 = n (2 g^\prime  -2) + \deg R\]
Let $\X_g$ be a genus $g\geq 2$ curve and $G$ its automorphism group (i.e, the group of automorphisms of the function field $\C (\X_g)$).  That   G is finite is shown in  \cite{nato-3} using Weierstrass points. 

Assume $|G|=n$.  Let $L$ be the fixed subfield of $\C (\X_g)$. The field extension $\C (\X_g)/L$ correspond to a finite morphism of curves $f : \X_g \to \X_g/G$ of degree $n$.   
Denote the genus of the quotient curve $\X_g/G$ by $g^\prime$ and $R$ the ramification divisor. Assume that the covering has $s$ branch points.  Each branch point $q$ has $n/e_P$ points in its fiber $f^{-1} (q)$, where $e_P$ is the ramification index of such points $P \in  f^{-1} (q)$.

Then, $R = \sum_{i=1}^s \frac n {e_P} \left( e_P - 1 \right)$. 
By the Riemann-Hurwitz formula we have 
\[  \frac  2 n \, (g-1) = 2 g^\prime - 2 + \frac 1 n \,  \deg R  = 2g^\prime - 2 + \sum_{i=1}^s \left( 1 - \frac 1 {e_P} \right) \] 
Since  $g \geq 2$ then the left hand side is $> 0$.  Then 
\[ 2 g^\prime - 2 +  \sum_{i=1}^s \left( 1 - \frac 1 {e_P} \right) \geq 0. \] The fact that $g^\prime$, $s$, and $e_P$ are non-negative integers implies that the minimum value of this expression is 1/42.  This implies that $n \leq 84 (g-1)$. 


Next we define another important invariant of the algebraic curves. Let $w_1, \dots , w_g$ be a basis of $H^0 (C, K)$ and $\gamma_1, \dots , \gamma_{2g}$ a basis for $H_1 (C, \z)$. The \textit{period matrix} $\Omega$ is the $g \times 2g$ matrix $ \Omega = \left[   \int_{\gamma_i} w_j \right]$.  The column vectors of the period matrix generate a lattice $\Lambda$ in $\C^g$, so that the quotient $\C^g/\Lambda$ is a complex  torus. This complex torus is called the \textit{Jacobian variety} of $C$ and denoted by the symbol $J (C)$. For more details see \cite{cornalba1} or \cite{fulton}.

If $\Lambda \subset \Lambda^\prime$ are lattices of rank $2n$ in $\C^n$ and $A= \C^n/\Lambda$, $B=\C^n/\Lambda^\prime$, then the induced map $A \to B$ is an \textit{isogeny} and $A$ and $B$ are said to be \textit{isogenous}.


The \textit{moduli space} $\M_g$ of curves of genus $g$ is the set of isomorphism classes of smooth, genus $g$ curves. $\M_g$ has a natural structure of a quasi-projective normal variety of dimension $3g-3$.  The compactification $\bar{\M_g}$ of $\M_g$ consists of isomorphism classes of stable curves.  A \textit{stable curve} is a curve whose only singularities are nodes and whose smooth rational components contain at least three singular points of the curve; see \cite[pg. 29]{cornalba1} and \cite[Chapter XII]{cornalba2}. $\bar{\M_g}$ is a projective variety 

Both $\M_g$ and $\bar{\M_g}$ are singular.  All the singularities arise from curves  with non-trivial automorphism group. It is precisely such curves that we intend to classify in this paper.

\subsection{Automorphism groups}\label{aut-groups}


\def\P{\mathbb P}

A curve $\X$ is called superelliptic if there exist an element $\tau \in \Aut( \X)$ such that $\tau$ is central and $g \left(\X / \< \tau \> \right) =0$.  Denote by $K$ the function field of $\X_g$ and assume that the affine equation of $\X_g$ is given some polynomial in terms of $x$ and $y$. 

Let $H=\< \tau \>$ be a cyclic subgroup of $G$ such that $| H | = n$ and $H \normal G$, where $n \geq 2$. Moreover, we assume that the quotient curve $\X_g / H$ has genus zero.   The \textbf{reduced automorphism group of $\X_g$ with respect to $H$} is called the group  $\G \, := \, G/H$, see \cite{super1}, \cite{Sa1}.  

Assume  $k(x)$ is the  genus zero subfield of $K$ fixed by $H$.   Hence, $[ K : k(x)]=n$. Then, the group  $\G$ is a subgroup of the group of automorphisms of a genus zero field.   Hence, $\G <  PGL_2(k)$ and $\G$ is finite. It is a classical result that every finite subgroup of $PGL_2 (k)$  is  isomorphic to one of the following: $C_m $, $ D_m $, $A_4$, $S_4$, $A_5$.

The group $\G$ acts on $k(x)$ via the natural way. The fixed field of this action is a genus 0 field, say $k(z)$. Thus, $z$
is a degree $|\G| := m$ rational function in $x$, say $z=\phi(x)$.   $G$ is a degree $n$ extension of $\G$ and $\G$ is a finite subgroup of $PGL_2(k)$.  Hence, if we know all the possible groups that occur as $\G$ then we can determine  $G$ and the equation for $K$. The list of all groups of superelliptic curves and their equations are determined in \cite{Sa1} and \cite{Sa2}.  

Let $C$ be a superelliptic curve given by the equation \[ y^n = f(x), \] where $\deg f = d$ and $\D (f, x) \neq 0$.  
Assume that $d> n$.  Then $C$ has genus \[ g = \frac 1 2 \left( n(d-1) -d - \gcd (n, d) \frac {} {}   \right) + 1\]
Moreover, if $n$ and $d$ are relatively prime then $g = \frac { (n-1) (d-1) } 2$,  see \cite{super1} for details.

\section{The majority of connected components of curves with non-trivial automorphisms correspond to superelliptic curves}

Let $g \geq 2$ be a given integer and $\X_g$ an irreducible, smooth, algebraic curve of genus $g$ defined over $\C$.  The automorphism group $\Aut(\X_g)$ is a finite group or order $\leq 84 (g-1)$.  For a given $g\geq 2$, the list of groups which occur as automorphism groups of genus $g$ curves can be determined following the methods in see \cite{kyoto} with some minor modifications. For hyperelliptic and superelliptic curves  we can determine no only the automorphisms groups (see \cite{issac}, \cite{Sa2}), but also their equations as in \cite{serdica}, \cite{Sa1}. 

The main question that we want to address in this section is: from the list of automorphism groups of curves for a fixed genus $g\geq 2$, how many of them come from superelliptic curves?   

The answer, at least for small genii, is that the majority of these groups come from superelliptic curves.  Hence, the majority of curves with non-trivial automorphism group have equation $y^n = f(x)$, for some $n$ and $f(x)$. Note that here the term "majority" means the majority of cases and not necessarily the "majority" in terms of dimension of these families.   

It is unclear if groups which come from non-superelliptic curves can be characterized as $g$ increases.  From the group theory point of view, these are groups $G$ of curves $\X_g$ such that for every  central element $\s \in G$, the genus of the quotient space   $g \left(\X_g / \< \s \> \right) \neq 0$. How do such groups look like when $g$ increases? 

Next we illustrate what happens for genus $g=4$.  The cases when $g=2, 3$ are easier and have appeared in the literature before.

\subsection{The case of genus 4}

In Table~\ref{tab_1}  we show all automorphism groups and their signatures for genus 4 algebraic curves. 
Each one of the families above is an algebraic locus in $\M_4$.   The data in this table was computed via the methods in \cite{kyoto}.

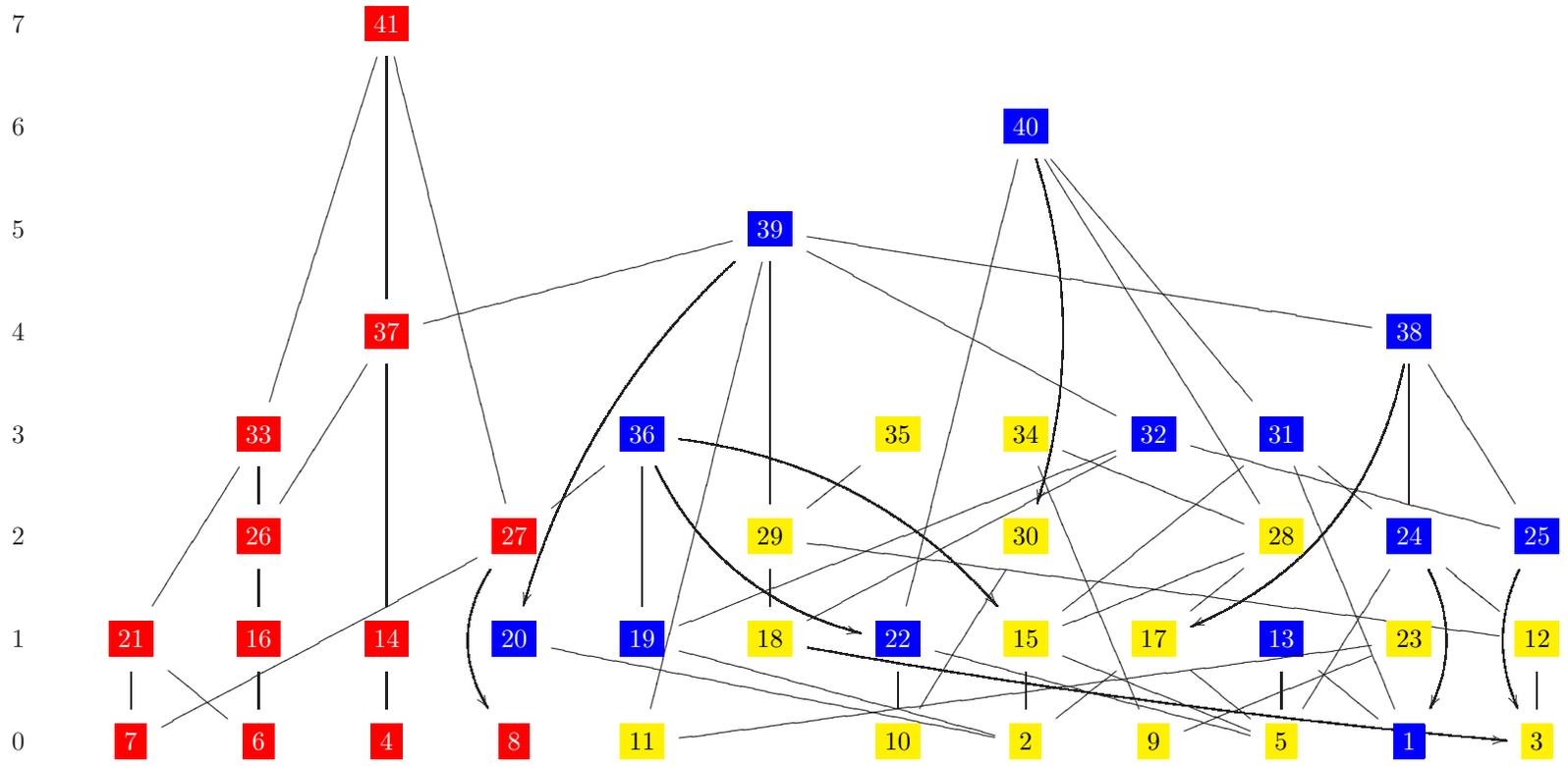
\begin{small}
\begin{table}[hbt] \label{tab_1}
\vskip -0.2cm
\begin{center}
\begin{tabular}{||c|c|c|c|c|c|c||}
\hline
\hline
$\#$ & dim &G& ID & sig & type & subs \\
\hline
\hline
1 & 0 & $S_5$ & (120,34) & 0-$(2, 4, 5)$ & 1 &  \\
2 & 0 &$C_3\times S_4$ &(72,42) & 0-$(2, 3, 12)$ & 3 &  \\
3 & 0 &         &(72,40) & 0-$(2, 4, 6)$ & 4 &  \\
4 & 0 & $V_{10}$ &(40,8) & 0-$(2, 4, 10)$ & 7 &  \\
5 & 0 & $C_6 \times S_3$ &(36,12) & 0-$(2, 6, 6)$ & 10 &  \\
6 & 0 & $U_8$    &(32,19) & 0-$(2, 4, 16)$ & 16 &  \\
7 & 0 & $SL_2(3)$ &(24,3) & 0-$(3, 4, 6)$ & 20 &  \\
8 & 0 & $C_{18}$  &(18,2) & 0-$(2, 9, 18)$ & 27 &  \\
9 & 0 & $C_{15}$ &(15,1) & 0-$(3, 5, 15)$ & 38 &  \\
10 & 0 & $C_{12}$  &(12,2) & 0-$(4, 6, 12)$ & 45 &  \\
11 & 0 & $C_{10}$  &(10,2) & 0-$(5, 10, 10)$ & 51 &  \\
12 & 1 & $S_3^2$&(36,10) & 0-$(2, 2, 2, 3)$ & 12 & 3 \\
13 & 1 & $S_4$&(24,12) & 0-$(2, 2, 2, 4)$ & 18 & 1, 2 \\
14 & 1 & $C_2\times D_5$ &(20,4) & 0-$(2, 2, 2, 5)$ & 21 & 4 \\
15 & 1 & $C_3\times S_3$ &(18,3) & 0-$(2, 2, 3, 3)$ & 30 & 2, 5 \\
16 & 1 & $D_8$           &(16,7) & 0-$(2, 2, 2, 8)$ & 35 & 6 \\
17 & 1 & $C_2\times C_6$ &(12,5) & 0-$(2, 2, 3, 6)$ & 46 & 2, 5 \\
18 & 1 & $C_2\times S_3$ &(12,4) & 0-$(2, 2, 3, 6)$ & 41 & 3 \\
19 & 1 & $A_4$       &(12,3) & 0-$(2, 3, 3, 3)$ & 43 & 2 \\
20 & 1 & $D_{10}$ &(10,1) & 0-$(2, 2, 5, 5)$ & 49 & 1 \\
21 & 1 & $Q_8$ &(8,4) & 0-$(2, 4, 4, 4)$ & 59 & 6, 7 \\
22 & 1 & $C_6$ &(6,2) & 0-$(2, 6, 6, 6)$ & 66 & 5, 10 \\
23 & 1 & $C_5$ &(5,1) & 0-$(5, 5, 5, 5)$ & 69 & 9, 11 \\
24 & 2 & $D_6$ &(12,4) & 0-$(2^{5})$ & 40 & 1, 5, 12 \\
25 & 2 & $D_4$ &(8,3) & 0-$(2^{4}, 4)$ & 57 & 3, 13 \\
26 & 2 & $D_4$ &(8,3) & 0-$(2^{4}, 4)$ & 56 & 4, 16 \\
27 & 2 & $C_6$ &(6,2) & 0-$(2^{3}, 3, 6)$ & 64 & 7, 8 \\
28 & 2 & $C_6$ &(6,2) & 0-$(2^{2}, 3^{3})$ & 65 & 15, 17 \\
29 & 2 & $S_3$ &(6,1) & 0-$(2^{2}, 3^{3})$ & 62 & 12, 18 \\
30 & 2 & $C_4$ &(4,1) & 0-$(2, 4^{4})$ & 77 & 10 \\
31 & 3 & $S_3$ &(6,1) & 0-$(2^{6})$ & 61 & 13, 15, 24 \\
32 & 3 & $V_4$ &(4,2) & 1-$(2, 2, 2)$ & 72 & 18, 19, 25 \\
33 & 3 & $C_4$ &(4,1) & 0-$(2^{4}, 4^{2})$ & 76 & 21, 26 \\
34 & 3 & $C_3$ &(3,1) & 0-$(3^{6})$ & 80 & 9, 28 \\
35 & 3 & $C_3$ &(3,1) & 0-$(3^{6})$ & 81 & 29 \\
36 & 3 & $C_3$  &(3,1) & 1-$(3, 3, 3)$ & 79 & 15, 19, 22, 27 \\
37 & 4 & $V_4$ &(4,2) & 0-$(2^{7})$ & 73 & 14, 26 \\
38 & 4 & $V_4$ &(4,2) & 0-$(2^{7})$ & 74 & 17, 24, 25 \\
39 & 5 & $C_2$ &(2,1) & 2-$(2, 2)$ & 82 & 11, 20, 29, 32, 37, 38 \\
40 & 6 & $C_2$ &(2,1) & 1-$(2^{6})$ & 83 & 22, 28, 30, 31, 38 \\
41 & 7 & $C_2$ &(2,1) & 0-$(2^{10})$ & 84 & 27, 33, 37 \\
\hline \hline
\end{tabular}
\end{center}
\caption{Hurwitz loci of genus 4 curves}
\end{table}
\end{small}

\clearpage




\begin{sideways}
\begin{minipage}{19.5cm}
\[
\xymatrixrowsep{.5cm}
\xymatrixcolsep{1cm}
\xymatrix@C-8pt{
 7        &    & & {\fcolorbox{red}{red}{\color{white}41}} \ar@{-}[ddddl] \ar@{-}[ddd] \ar@{-}[dddddr]&    &    &    &    &    &    &    &              \\
 6        &    &    &    &    &    &    &    &      {\fcolorbox{blue}{blue}{\color{white}40}} \ar@{-}[dddddl] \ar@/^1.2pc/[dddd] \ar@{-}[dddrr] \ar@{-}[ddddrr]&    &    &      \\
 5        &    &    &    &    &    &
 {\fcolorbox{blue}{blue}{\color{white}39}} \ar@{-}[dddddl]   \ar@{-}[dlll] \ar@{-}[ddd] \ar@/_1.3pc/[ddddll] \ar@{-}[ddrrr] \ar@{-}[drrrrr]&    &    &    &    &     \\
 4        &    &    &
 {\fcolorbox{red}{red}{\color{white}37}} \ar@{-}[ddl]\ar@{-}[ddd] &     &    &    &    &    &    &    & {\fcolorbox{blue}{blue}{\color{white}38}}  \ar@{-}[dd]\ar@{-}[ddr] \ar@/^2.0pc/[dddll]  \\
 3        &    &
 {\fcolorbox{red}{red}{\color{white}33}} \ar@{-}[d] \ar@{-}[ddl] &    &    &    {\fcolorbox{blue}{blue}{\color{white}36}} \ar@{-}[dl] \ar@{-}[dd]\ar@/_1.5pc/[ddrr] \ar@/^1.6pc/[ddrrr]&    & {\fcolorbox{yellow}{yellow}{\color{black}35}}  \ar@{-}[dl] & {\fcolorbox{yellow}{yellow}{\color{black}34}}   \ar@{-}[drr] \ar@{-}[dddr]& {\fcolorbox{blue}{blue}{\color{white}32}}  \ar@{-}[drrr] \ar@{-}[ddllll]  \ar@{-}[ddlll] & {\fcolorbox{blue}{blue}{\color{white}31}}  \ar@{-}[dddr] \ar@{-}[dr] \ar@{-}[ddll]&      \\
2   &    &
  {\fcolorbox{red}{red}{\color{white}26}} \ar@{-}[d] &    &
     {\fcolorbox{red}{red}{\color{white}27}} \ar@/_1.5pc/[dd] \ar@{-}[ddlll]    &    &
     {\fcolorbox{yellow}{yellow}{\color{black}29}}  \ar@{-}[d]\ar@{-}[drrrrrr] &    & {\fcolorbox{yellow}{yellow}{\color{black}30}}  \ar@{-}[ddl] &    &
 {\fcolorbox{yellow}{yellow}{\color{black}28}}  \ar@{-}[dl] \ar@{-}[dll]& {\fcolorbox{blue}{blue}{\color{white}24}}  \ar@/^1.2pc/[dd] \ar@{-}[ddl] \ar@{-}[dr]& {\fcolorbox{blue}{blue}{\color{white}25}}  \ar@/_1.1pc/[dd]       \\
1        &
 {\fcolorbox{red}{red}{\color{white}21}} \ar@{-}[d] \ar@{-}[dr]&
 {\fcolorbox{red}{red}{\color{white}16}} \ar@{-}[d] &
 {\fcolorbox{red}{red}{\color{white}14}} \ar@{-}[d] &
 {\fcolorbox{blue}{blue}{\color{white}20}} \ar@{-}[drrrr] &
 {\fcolorbox{blue}{blue}{\color{white}19}} \ar@{-}[drrr] &
 {\fcolorbox{yellow}{yellow}{\color{black}18}} \ar@/_0.3pc/[drrrrrr] &
 {\fcolorbox{blue}{blue}{\color{white}22}} \ar@{-}[d]\ar@{-}[drrr] &
 {\fcolorbox{yellow}{yellow}{\color{black}15}} \ar@{-}[d]\ar@{-}[drr] &
 {\fcolorbox{yellow}{yellow}{\color{black}17}} \ar@{-}[dr]\ar@{-}[dl] &
 {\fcolorbox{blue}{blue}{\color{white}13}} \ar@{-}[d]\ar@{-}[dr] &
 {\fcolorbox{yellow}{yellow}{\color{black}23}} \ar@{-}[dllllll] \ar@{-}[dll]  &
 {\fcolorbox{yellow}{yellow}{\color{black}12}} \ar@{-}[d]&     \\
0    &
{\fcolorbox{red}{red}{\color{white}7}}  &
{\fcolorbox{red}{red}{\color{white}6}}  &
{\fcolorbox{red}{red}{\color{white}4}}  &
{\fcolorbox{red}{red}{\color{white}8}}    &
{\fcolorbox{yellow}{yellow}{\color{black}11}}  &    &
{\fcolorbox{yellow}{yellow}{\color{black}10}}      &
{\fcolorbox{yellow}{yellow}{\color{black}2}}   &
{\fcolorbox{yellow}{yellow}{\color{black}9}}    &
{\fcolorbox{yellow}{yellow}{\color{black}5}}   &
{\fcolorbox{blue}{blue}{\color{white}1}}   &
{\fcolorbox{yellow}{yellow}{\color{black}3}}      \\
}
\]
\end{minipage}
\end{sideways}

\smallskip

In the  diagram above the red and yellow entries denote the superelliptic curves (hyperelliptic and non-hyperelliptic respectively).  Notice that from 41 total cases, only 13 are non-hyperelliptic.

\def\l{\lambda }


\begin{table}[hbt]
\caption{Equations of genus 4 superelliptic curves}\label{tab_eq}
\begin{center}
\begin{tabular}{||c|c|c|c||}
\hline \hline
$\#$ & dim & aut & equation \\
\hline \hline
&&&\\
\textbf{23} & 1 & (5,1) &  $y^5=x(x-1)(x-\l)$ \\
9 & 0 & (15,1) &  $y^5=x^3-1$  \\
11 & 0 & (10,2) &  $y^5=x(x^2-1)$  \\  
&&&\\
\hline  &&&\\
\textbf{34} & 3 & (3,1) &  $y^3= x(x-1)(x-\a_1)(x-\a_2)(x-\a_3)$ \\   
28 & 2 & (6,2) &  $ y^3= (x^2-1)(x^2-\a_1)(x^2-\a_2)$ \\  
15 & 1 & (18,3) & $y^3 = x^6+ \l x^3 +1$ \\
17 & 1 & (12,5) &  $y^3= (x^2-1) (x^4 - \lambda x^2 +1 )$ \\
2 & 0 & (72,42) & $y^3 = x (x^4-1)$   \\
5 & 0 & (36,12) & $y^3 = x^6-1 $   \\ 
&&&\\
\textbf{35} & 3 & (3,1) &  $y^3=(x^2-2)(x^4+bx^2+cx+d)$ \\  
29 & 2 & (6,1) &  $y^3-1= x \left( x^5+(b-2)x^3+x^3c-(2b+1/2)x-2c \right)$ \\ 
12 & 1 & (36,10) &  $y^3-1=x^6+\l x^3 +1$  \\
18 & 1 & (12,4) & $y^3-1=(x^2-1)(x^2-\a_1)(x^2-\a_2)$  \\
3 & 0 & (72,40) &   $y^3-1  =  x^6- 1$    \\ 
&&&\\
\textbf{22} & 1 & (6,2) &  $y^6= x(x-1)(x-\a)$ \\ 
30 & 2 & (4,1) &  $y^4=x^2(x-1)(x-\a_1)(x-\a_2)$  \\
10 & 0 & (12,2) & $y^4= x^2 (x^3-1)$    \\  
&&&\\
\hline  &&&\\
\textbf{41} & 7 & (2,1) &  $y^2=f(x), \, \, \deg f = 9, 10$ \\
37 & 4 & (4,2) & $ y^2= x^{10} + a_1 x^8 + a_2 x^6 + a_3 x^4 + a_4 x^2 +1$ \\
33 & 3 & (4,1) & $ y^2 = x (x^8 + a_1 x^6 + a_2 x^4 + a_3 x^2 +1)$ \\
26 & 2 & (8,3) & $ y^2 =x (x^4+\l_1 x^2 +1)(x^4+\l_2 x^2 +1) $ \\
27 & 2 & (6,2) & $ y^2 =x^9+a_1x^6+a_2 x^3 +1 $ \\
4 & 0 & (40,8) & $ y^2 =x^{10}-1 $  \\
6 & 0 & (32,19) & $ y^2 =x (x^8-1) $  \\
7 & 0 & (24,3) & $ y^2 = x (x^4-1)(x^4+2i \sqrt{3} x^2 +1) $ \\
8 & 0 & (18,2) & $ y^2 =x^9+1 $  \\
21 & 1 & (8,4) & $ y^2 =x (x^4-1)(x^4 + \l x^2 +1) $ \\
14 & 1 & (20,4) & $ y^2 =x^{10} + \l x^5 +1 $ \\
16 & 1 & (16,7) & $ y^2 =x(x^8 + \l x^4 +1) $ \\
&&&\\
\hline \hline
\end{tabular}
\end{center}
\end{table}

For each of the superelliptic cases we can write explicitly the equation of each family of curves. The equations in all other cases are not known.

As seen from the tables, the superelliptic curves are well understood for genus $g=4$ even though many questions remain.  Moreover, we could  compile such tables for any $g>2$, even though as $g$ increases the computations become more challenging. 

Given a  curve of genus 4, can we determine if the curve  belongs to any of the cases of the above table?   Moreover, for any of the curves in the above families defined over $\overline Q$ can we determine if such curve is isomorphic to a curve defined over $\Q$?  Some of these questions will be treated in the next section. 

\section{Equations of superelliptic curves over their minimal field of definition}

For each $g$, the moduli space   $\M_g$  (resp., $\H_g$)  is the set of isomorphism classes of genus  $g$
algebraic   (resp.,  hyperelliptic) curves  $\X_g$ defined over an algebraically closed field  $k$. It is well
known that $\M_g$ (resp.,  $\H_g$) is a $3g-3 $  (resp., $2g-1$) dimensional variety.   Let $L$ be a subfield of
$k$. If $\X_g$ is  a genus $g$ curve defined  over $L$, then clearly $[\X_g]\in  \M_g (L)$. Generally, the converse
does not hold. In  other words,  the moduli spaces  $\M_g$ and $\H_g$ are coarse moduli spaces.

Let $\X$ be a curve defined over $k$.  A field $F \subset k$  is called a \textbf{  field of definition} of $\X$
if there exists $\X'$ defined over $F$ such that $\X'$ is isomorphic to $\X$ over $k$.

The \textbf{ field of moduli} of $\X$ is a subfield $F \subset k$ such that for every automorphism $\sigma$ of $k$
$\X$ is isomorphic to $\X^\sigma$  if and only if $\, \,  \sigma_F = id$.

We will use $\p=[\X]\in   \M_g$ to denote  the corresponding \textbf{ moduli point} and $\M_g (\p)$ the
residue field of $\p$  in $\M_g$. In characteristic zero the field of moduli of $\X$ coincide with the
residue field $\M_g (\p)$ of the point $\p$  in $\M_g$; see Baily \cite{Ba2}. 
 The notation  $\M_g (\p)$ (resp., $M(\X)$ ) will be used to denote
the field of moduli of $\p \in \M_g$ (resp., $\X$). If there is a curve $\X^\prime$ isomorphic to $\X$ and
defined over $M(X)$, we say that $\X$ has a \textbf{rational model over its field of  moduli}.

As mentioned above, the field of moduli of curves is not necessarily a field of definition. However, finding such
curves is not easy. Shimura was the first to provide an example of a family of hyperelliptic curves such that the
field of moduli is not a field of definition; see \cite{Shi}.

\begin{exa} [Shimura]
The family of curves is given by
\[Y^2= a_0 X^m + \sum_{r=1}^m (a_r X^{m+r} + (-1)^r {\bar a}_r X^{m-r})\]
where $m$ is odd, $a_m=1, a_0\in \mathbb R$,  $a_1, \dots , a_{m-1} $ are complex numbers and ${\bar a}_r$ denotes
the complex conjugation of $a_r$. These curves have complex conjugation, so the field of moduli is at most $\R$.
However, for an appropriate choice of $a_i$'s ( algebraically independent over $\Q$) these curves are not defined
over $\R$.
\end{exa}

It is well known that for a given algebraic curve $\X$ defined over any subfield of $\, \C$, there exists a
field of definition for $\X$ which is a finite algebraic extension of the moduli field  $M(\X) $. We
naturally are interested in curves for which the field of moduli is a number field. Then, their field of
definition can be also chosen to be a number field. For the proof of the following result see \cite{Wo},
Theorem 3.
\begin{thm}
A compact Riemann surface is a Belyi surface if and only if its moduli field is a number field.
\end{thm}

What are necessary conditions for a curve to have a rational model over its field of moduli? We consider only
curves of genus $g> 1$; curves of genus 0 and 1 are known to have a rational model over its field of moduli. In
(1954) Weil  showed that;

\begin{itemize}
\item [i)]  {\it For every  curve $\X$ with trivial automorphism group  the field of moduli is a field of
definition.}
\end{itemize}

Later work of  Belyi, Shimura, Coombes-Harbater, D\'ebes, Douai, Wolfart et al. has added other conditions which
briefly are summarized below.

\noindent The field of moduli of a curve $\X$ is a field of definition if:

\begin{itemize}
\item [  ii)]   {\it $Aut (\X)$ has no center and has a complement in the automorphism group of $Aut (\X)$}

\item [  iii)]  {\it The field of moduli $M (\X)$ is of cohomological dimension $\leq 1$}

\item [  iv)]  {\it The canonical $M ( \X)$-model of $\X / Aut(\X)$ has $M (\X)$-rational points.}
\end{itemize}


\def\bP{\mathbb P}
\def\V{\mathbf V_4}
\def\bAut{  {\overline{\mbox {Aut}} \, \,}}
\def\<{\langle}
\def\>{\rangle}
\def\iso{{\, \cong\, }}
\def\S{\mathcal S}

\subsection{Field of moduli of superelliptic curves}

It seems that superelliptic curves are  the most interesting examples on the field of moduli problem. Shimura's
family and Earle's family of curves (i.e., with non-trivial obstruction) are both families of hyperelliptic
curves.   Moreover, the following is proved in \cite{Shi}.

\begin{thm}[Shimura] No generic hyperelliptic curve of even genus has a model rational over its field of
moduli.
\end{thm}

\noindent Consider the following problem:

\begin{prob}  Let the  moduli point $\p \in \S_g$ be given, where $S_g$ is one of the loci of the superelliptic curves in the corresponding moduli space $\M_g$. Find necessary and sufficient conditions
that the field of moduli $\S_g (\p)$ is a field of  definition.  If $\p$ has a rational model $\X_g$ over its
field of moduli, then determine explicitly the equation of $\X_g$.  
\end{prob}

In 1993, Mestre solved the above problem for genus two curves with automorphism group of order 2. He showed that
given a   field $L \subset \C$, there is a bijective correspondence between the $L$-rational points of $\M_2$ and
$\bar L$-isomorphism classes of pairs $(M,D)$ ($M$ being a genus zero curve over $L$ and $D$ an $L$-rational
effective divisor of degree $6$ without multiplicities). Moreover, under this correspondence, a point $\p \in
\M_2$ such that $|Aut(\p)|=2$  is given by a curve $C$ defined over $L$ if and only if $M$ is isomorphic to
$\bP\sp 1 \sb L$. The proof is based on the results of Clebsch, Bolza, and Igusa  on the classical theory of
invariants of binary sextics.

In \cite{Sh3}, following a different approach,  we show that for genus 2 curves with extra automorphisms   the
field of moduli is a field of definition. 

Following the same approach as Mestre and using  a coordinate for $\H_3$ in terms of absolute invariants $t_1, \dots , t_6$, as defined in \cite{mod-g-3} one would like to have an algorithm which determines when the field of moduli is a field of definition for the generic hyperelliptic curve of genus 3.  Moreover,   an equation of the curve over its minimal field of definition would be desirable.

\begin{prob}
Let $\p = (t_1, \dots , t_6)$ be a point in the genus 3 hyperelliptic moduli, where $t_1, \dots , t_6$ are defined as in \cite{mod-g-3}.  Provide an explicit  equation of a curve in terms of $t_1, \dots , t_6$,  when possible. 
\end{prob}

\subsection{Generic superelliptic curves}  Invariants of binary forms are known up to binary decimics and theoretically we understand the general case ($n > 5$).  Such invariants give a description of a point in the corresponding moduli space and should be enough to determine such equations for "generic" superelliptic curves of "small" genus.  The natural question is, how far can be pushed the computational limits and up to what genus $g$ do we get explicit results?
 Superelliptic curves with extra automorphisms are curves $y^n=f(x)$ with automorphism group of order $> n$.  Such curves have to be treated differently from the generic case. Probably techniques used for genus 2 and 3 should work for higher genii as well.

\section{Curves with minimal height}

Let us assume that from above section we have determined a Weierstrass equation of some superelliptic curve defined over 
some number field $K$.  Let  $\O_K$ be its ring of integers. The main question is how "good" is this equation?  

In general, the algorithms that we have from section 2 provide equations with extremely large coefficients.  Can we fid an equation which is "minimal" in some sense. Below we give briefly some preliminaries on heights of curves and the moduli height of a given curve. For more details please see \cite{nato-12}. 

Let $\X_g$ be an irreducible algebraic curve with affine equation $F(x, y)=0$ for $F(x, y) \in K [x, y]$.  We define
the \textbf{height of the  curve over $K$} to be
\[H_K(\X_g):= \min \left \{    H_K(G) \, : H_K(G) \leq H_K(F) \right \}. \]
where the curve $G(x, y) =0$ is isomorphic to $\X_g$ over $K$. 

If we consider the equivalence over $\bar K$ then we get another height which we denote it as $\overline H_K (\X_g)$ and call it \textbf{the height over the algebraic closure}. Namely, 
\[ \overline H_K(\X_g)= \min \{H_K(G): H_K(G) \leq H_K(F)\},\]
 where the curve $G(x, y)=0$ is isomorphic to $\X_g$ over $\overline K$.
 
In the case that $K=\Q$ we do not write the subscript $K$ and use $H(\X_g)$ or  $\overline H(\X_g)$.  Obviously, for any algebraic curve $\X_g$ we have $\overline H_K(\X_g) \leq H_K(\X_g)$.  For the proof of the following results see \cite{nato-8}.

\begin{lem} Let $K$ be a number field such that $[K:\Q] = d$.  Then, 
$H_K(\X_g)$  and   $\overline H_K(\X_g)$    are  well defined.
\end{lem}

\begin{thm}
Let $K$ be a number field such that $[K:\Q] \leq d$. Given a constant $c$  there are only finitely many curves (up to isomorphism)  such that $H_K(\X_g) \leq c$.
\end{thm}

\subsection{Moduli height of curves}

Let $g$ be an integer $g \geq 2$ and $\M_g$ denote the coarse moduli space of smooth, irreducible algebraic curves of genus $g$. It is known that $\M_g$ is a quasi projective variety of dimension $3g-3$.  Hence, $\M_g$ is embedded in $\P^{3g-2}$. Let $\p \in \M_g$. We call the moduli height $\mH(\p)$ the usual height $H(P)$ in the projective space $\P^{3g-2}$.  Obviously, $\mH(\p)$ is an invariant of the curve.

\begin{thm}For any constant $c\geq 1$, degree $d\geq 1$, and genus $g\geq 2$  there are finitely many superelliptic curves $\X_g$ defined over the ring of integers $\O_K$ of an algebraic number field $K$ such that   $[K:\Q] \leq d$ and  $\mH (\X_g) \leq c$.
\end{thm}

\subsection{Classifying superelliptic curves}

One of the main goals would be to classify all superelliptic curves over $\Z$ for a fixed genus $g$.  Of course the same question can be asked for any rink of integers $\O_K$. 

One can attempt to write down list of such curves as the tables in \cite{nato-8} for genus 2 curves of height 1.  However, what is missing in those tables is that for every moduli point $\p \in \M_2$ is given a curve of minimal height corresponding to $\p$ but not all its twists.  In other words, there are usually other curves defined over $\Z$ belonging to the same moduli point.  In \cite{bin} it is intended a method of how to list all such curves for every moduli point.

Another way of classifying or listing the curves would be through their moduli height.  This would make more sense theoretically because the height of the moduli point is a more natural invariant.  However, computationally this definitely would be more challenging.  Curves with small moduli height don't necessarily have small height (i.e., small coefficients).

\nocite{*}

\bibliographystyle{} 

\bibliography{mybib}{}


\end{document}